\documentclass[a4paper,12pt,reqno,twoside]{amsart}

\usepackage[utf8]{inputenc} 		
\usepackage[T1]{fontenc} 
\usepackage{mathpazo}
\usepackage{amssymb}
\usepackage{amsmath}
\usepackage{graphicx}
\usepackage[colorlinks,citecolor=blue,urlcolor=blue]{hyperref}
\usepackage{cite}
\usepackage[hmargin=2.5cm,vmargin=2.5cm]{geometry}
\usepackage{mathrsfs} 

\newenvironment{myproof}{\paragraph{Proof:}}{\hfill$\blacksquare$}

\newcommand{\be}{\begin{eqnarray}}
\newcommand{\ee}{\end{eqnarray}}

\newcommand{\beq}{\begin{equation}}
\newcommand{\eeq}{\end{equation}}
\newcommand{\beqn}{\begin{equation*}}
\newcommand{\eeqn}{\end{equation*}}

\newcommand{\round}[1]{\lfloor#1\rfloor}

\DeclareMathOperator{\tr}{tr}

\newtheorem{thm}{Theorem}[section]

\newtheorem{cor}[thm]{Corollary}
\newtheorem{lem}[thm]{Lemma}

\newtheorem{fact}[thm]{Fact}

\newcommand\cB{{\mathcal B}}
\newcommand\cC{{\mathcal C}}

\newcommand\cL{{\mathcal L}}

\newcommand\cN{{\mathcal N}}

\newcommand\cP{{\mathcal P}}

\newcommand\bA{{\mathbb A}}

\newcommand\bL{{\mathbb L}}

\newcommand\bN{{\mathbb N}}

\newcommand\bR{{\mathbb R}}

\newcommand{\ve}{\varepsilon}

\newtheorem{claim}{Claim}

\renewcommand{\qed}{\hfill\blacksquare}

\begin{document}


\title[functional correlation decay for non-uniformly expanding maps]{functional correlation decay and multivariate normal approximation for non-uniformly expanding maps}

\author[Juho Lepp\"anen]{Juho Lepp\"anen}
\address[Juho Lepp\"anen]{
Department of Mathematics and Statistics, P.O.\ Box 68, Fin-00014 University of Helsinki, Finland.}
\email{juho.leppanen@helsinki.fi}

\maketitle

 

\begin{abstract}
In the setting of intermittent Pomeau-Manneville maps with time dependent parameters, we show a functional correlation bound widely useful for the analysis of the statistical properties of the model. We give two applications of this result, by showing that in a suitable range of parameters the bound implies the conditions of the normal approximation methods of Stein and Rio. For a single Pomeau-Manneville map belonging to this parameter range, both methods then yield a multivariate central limit theorem with a rate of convergence.
\end{abstract}

\section{Introduction} 

Recently, general methods have been devised for obtaining rates of convergence in the multivariate CLT. Given a measure preserving transformation \(T : \, X \to X\) on a probability space \((X,\cB,\mu) \) and a function \(f : \, X \to \bR^d\), \(d \ge 1\), with \(\int f \, d\mu = 0 \), the sequence \((f \circ T^n)_{n \ge 1}\) is a centered stationary process. We say that it satisfies the central limit theorem (CLT), if the normalized Birkhoff sums \(N^{-1/2} \sum_{k=0}^{N-1} f \circ T^k \) converge in distribution to a \(d\)-dimensional Gaussian random variable. In \cite{hella2016}, certain correlation-decay conditions were formulated by using Stein's method of normal approximation. Once the system is shown to satisfy these conditions, multivariate CLT augmented by a rate of convergence follows immediately. By building on a method due to Rio \cite{rio1996}, different conditions for the same purpose were established by F. P\`{e}ne in \cite{pene2005}.\\
\indent In both of the above papers, the authors verify the formulated conditions for a system with exponential decay of correlations, namely the Sinai billiard. In \cite{pene2005}, also the Knudsen gas was considered. A natural question arising is whether these conditions are satisfied also by systems that exhibit a weaker rate of correlation decay, say polynomial. The purpose of this paper is to show that the answer is positive in the case of the intermittent Pomeau-Manneville maps of a suitable parameter range. We demonstrate how these maps satisfy a certain functional correlation bound (Theorem \ref{corr}) from which the conditions of \cite{pene2005,hella2016} readily follow. As results, we obtain two versions of the multivariate CLT with speed in this setting of Pomeau-Manneville maps.   

\subsection{Some notation.} Given a probability space \((X,\cB,\mu)\) and a function \(f: \, X \to \bR^d\), we denote \(\mu(f) = \int_{X} f \, d\mu\). The Lebesgue measure on the unit interval is denoted by \(m\). The coordinate functions of \(f\) are denoted by \(f_{\alpha}\), $\alpha \in\{1,\dots,d\}$, and we let
\beqn
\|f\|_\infty = \max_{1\le\alpha\le d}\|f_\alpha\|_\infty.
\eeqn

\indent We endow \(\bR^d\) with the max-norm \(\Vert x \Vert_{\infty} = \max_{\alpha =1,\ldots ,d} \vert x_{\alpha} \vert \), and for a Lipschitz continuous function \(f \, : [0,1] \to \bR^{d}\) define 
\begin{align*}
\text{Lip}(f) =  \max_{\alpha = 1,\ldots,d}\sup _{x \neq y} \frac{ | f_{\alpha}(x) - f_{\alpha}(y) | }{| x-y |},
\end{align*}
and \(\Vert f \Vert_{\text{Lip}} = \Vert f \Vert_{\infty} + \text{Lip}(f)\). \\
\indent Given a function \(F: \, [0,1]^k \to \bR\) and \(i \in \{1,\ldots, k\}\), we denote by \(\text{Lip}(F;i)\) the quantity
\begin{align*}
\sup_{y_1,\ldots, y_d \in [0,1]}\sup_{a_i \neq b_i} \frac{|F(y_1,\ldots, y_{i-1}, a_i, y_{i+1},\ldots , y_d) - F(y_1,\ldots, y_{i-1}, b_i, y_{i+1},\ldots , y_d)|}{|a_i-b_i|},
\end{align*}
 and say that \(F\) is Lipschitz continuous in the \(i\)th coordinate \(x_i\), if \(\text{Lip}(F;i) < \infty\).
\subsection{A class of intermittent maps}\label{intermaps} Following \cite{liverani1999}, for each \(\alpha \in (0,1)\), we define the Pomeau-Manneville map \(T_{\alpha} : \, [0,1] \to [0,1] \) by
\begin{align*}
T_{\alpha }(x) = \begin{cases} x(1+ 2^{\alpha }x^{\alpha}) & \forall x \in [0, 1/2), \\
2x-1 & \forall x \in [1/2,1].
 \end{cases}
\end{align*}

The fundamental characteristic of these maps is the intermittent behaviour they exhibit due to the presence of the neutral fixed point at the origin, where their derivative equals one. At every other point, the maps expand locally uniformly. Mappings with such properties are known to demonstrate a polynomial rate of correlation decay \cite{liverani1999,hu2004,young1999}. Moreover, this rate is known to be sharp \cite{gouezel2004}. It follows from \cite{liverani1999}, that each \(T_{\alpha}\) admits an invariant SRB measure \(\hat{\mu}_{\alpha}\), whose density belongs to the convex cone \\
\beqn
\begin{split}
\cC_*(\alpha) = \{f\in C((0,1])\cap L^1\,:\, & \text{$f\ge 0$, $f$ decreasing,} 
\\
& \text{$x^{\alpha+1}f$ increasing, $f(x)\le 2^{\alpha} (2 + \alpha) x^{-\alpha} m(f)$}\}.
\end{split}
\eeqn

Limit theorems for Pomeau-Manneville maps have been established in several earlier papers; see, for instance, \cite{nicol2016,dedecker2009,pollicott2002,zweinmuller2003,bahsoun2014,gouezel2007}. Results on the rate of convergence in the univariate CLT were obtained by S. Gou{\"e}zel in \cite{gouezel2005} by implementing a general Young-tower technique. For parameters \(\alpha < 1/3\), Gou{\"e}zel's results show that for any Hölder continuous function \(f : [0,1] \to \bR\) with \(\int f \, \hat{\mu}_{\alpha} = 0\), the scaled time average \(N^{-1/2} \sum_{k=0}^{N-1} f \circ T_{\alpha}^k \) converges in distribution at the optimal rate \(N^{-1/2}\) to normal distribution \(\cN(0,\sigma^2)\) with mean zero and variance \(\sigma^2 > 0\). Even for \(\alpha \in [1/3, 1/2)\), Gou{\"e}zel establishes a polynomial speed in the CLT when some additional control is given on the behavior of \(f\) around the origin. If \(\alpha \ge 1/2\), then it is seen from \cite{gouezel2004_2} that the limit still exists, and the limiting distribution is, depending on the properties of \(f\), either a normal distribution or a stable distribution.\\
\indent In the present paper we consider a Pomeau-Manneville map \(T_{\alpha}\) with \(\alpha < 1/3\), and establish a rate of convergence in the CLT for multivariate functions \(f\). We emphasize that although multivariate CLTs follow from one-dimensional CLTs by a well-known general argument, convergence rates obtained in dimension one do not directly transfer to higher dimensions.

\subsection{Functional correlation decay.}  We fix once and for all a real number \(\beta_* \in (0,1)\), and call a sequence of mappings \((T_{\alpha_n})_{n\ge1}\) admissible, if \(\alpha_n \le \beta_*\) for all \(n \ge 1\). Given such a sequence of parameters \((\alpha_n)\),  we abbreviate \(\cC_* = \cC_*(\beta_*)\), \(T_n = T_{\alpha_n}\) and \( \widetilde{T}_{n} = T_n \circ \cdots \circ T_1\).  \\
\indent For comparing quantities, we introduce the following notations. For any real-valued functions \(f\) and \(g\), we denote \(g(x) \lesssim_{\theta} f(x)\) if there exists \(C > 0\) depending only on \(\theta\) with \(g(x) \le Cf(x)\). Moreover, \(\lesssim\) means \(\lesssim_{\beta_*}\) and \(g(x) \sim_{\theta} f(x)\) means that \(g(x) \lesssim_{\theta} f(x)\) and \(f(x) \lesssim_{\theta} g(x)\). \\
\indent Here is the main result:

\begin{thm}\label{corr} Let \((T_n)_{n\ge 1}\) be any admissible sequence of mappings. Let \(F: \, [0,1]^{k+1} \to \bR\) be a bounded function, and fix integers \(0 = n_0 \le n_1 \le \ldots \le n_k\), \(1 \le l_1 < \ldots < l_p < k\). Suppose that \(F\) is Lipschitz continuous in the coordinate \(x_{\alpha}\) whenever \(1 \le \alpha \le l_p +1\),  and denote by \(H(x_0,\ldots, x_{p})\) the function
\begin{align*}
F(x_0,\widetilde{T}_{n_1}(x_0),\ldots, \widetilde{T}_{n_{l_1}}(x_0), \widetilde{T}_{n_{{l_1}+1}}(x_1),\ldots,\widetilde{T}_{n_{{l_2}}}(x_1), \ldots, \widetilde{T}_{n_{l_p+1}}(x_{p}),\ldots, \widetilde{T}_{n_k}(x_{p})).
\end{align*}
Then, for any probability measures \(\mu,\mu_1,\ldots,\mu_p \) whose densities belong to \(\cC_*\),
\begin{align}
&\left| \int H(x,\ldots,x) \, d\mu(x) - \idotsint H(x_0,\ldots,x_{p}) \, d\mu(x_0)d\mu_1(x_1)\ldots \, d\mu_p(x_{p}) \right| \label{coreq}\\ 
&\lesssim (\Vert F \Vert_{\infty} + \max_{1 \le\alpha \le l_p+1} \textnormal{Lip}(F;\alpha))\sum_{i=1}^p \rho(n_{l_{i}+1}-n_{l_i}), \notag
\end{align}
where \(\rho(n) = n^{-\frac{1}{\beta_*}+1}(\log n)^{\frac{1}{\beta_*}}\) for \(n \ge 2\), and \(\rho(0) = \rho(1) = 1\).
\end{thm}

In the special case \(p=1\), the function \(H\) in Theorem \ref{corr} becomes
\begin{align*}
H(x,y) = F(x,\widetilde{T}_{n_1}(x),\ldots, \widetilde{T}_{n_{l}}(x), \widetilde{T}_{n_{{l}+1}}(y),\ldots,\widetilde{T}_{n_{k}}(y)),
\end{align*}
where \(l = l_1\). Informally speaking, the result then states that 
\begin{align*}
&\int H(x,x) \, d\mu(x) \approx \iint H(x,y) \, d\mu(x) \, d\mu_1(y),
\end{align*}
where the error in this approximation is polynomial in the gap \(n_{l+1} - n_l\) between \(n_l\) and \(n_{l+1}\). The general formulation concerns the case where there are \(p\) gaps \(n_{i+1} - n_i\), located at \(i = l_1,\ldots, l_p\). In this case too, we obtain a result stating that the integral \(\int H(x,\ldots, x) \, d\mu(x)\) almost factorizes to a product of one-dimensional integrals. 

Our motivation for proving Theorem \ref{corr} lies in its implications to normal approximation of scaled Birkhoff sums. Indeed, the applicability of the general normal approximation methods due to Stein \cite{hella2016} and Rio \cite{pene2005} completely depend on the ability to control quantities of the form \eqref{coreq} for various functions \(F\).  The general form of \(F\) in the theorem enables us to verify the correlation-decay conditions of these methods with minimal effort, indicating that Theorem \ref{corr} is a rather versatile tool for proving limit theorems in the setting of Pomeau-Manneville maps. These applications concern a single map \(T_{\beta_*}\) only, but we have decided to formulate and prove the result for a sequences of maps \((T_n)_{n\ge 1}\) instead, for this more general form will be useful to showing limit theorems beyond this paper. One example of such a limit theorem is the following multicorrelation bound, which we obtain immediately from Theorem \ref{corr} by taking \(F\) to be a product of one-dimensional observables.

\begin{cor}\label{multicor} Let \((T_n)_{n \ge 1}\) be an admissible sequence of mappings, \(f_0,\ldots, f_{l}\) Lipschitz continuous functions \([0,1] \to \bR\), and \(f_{l+1},\ldots,f_k \in L^{\infty}([0,1])\). Fix integers \(0 = n_0 \le n_1 \le \ldots \le n_k\), and denote
\begin{align*}
H &= f_0 \cdot f_1 \circ \widetilde{T}_{n_1} \cdots f_{l}\circ \widetilde{T}_{n_{l}} \\
G &=  f_{l+1} \circ \widetilde{T}_{n_{{l}+1}} \cdots f_{k}\circ \widetilde{T}_{n_{k}}.
\end{align*}
Then, for any probability measure \(\mu\) with density \(h \in \cC_*\),
\begin{align}
&\left| \int HG \, d\mu - \int H \, d\mu \int G \, d\mu \right| \notag\\ 
&\lesssim  \prod_{i=0}^l \Vert f_i \Vert_{\textnormal{Lip}}\prod_{i=l+1}^k \Vert f_i \Vert_{\infty} \rho(n_{l+1}-n_{l}), \label{eq:cor}
\end{align}
where \(\rho(n) = n^{-\frac{1}{\beta_*}+1}(\log n)^{\frac{1}{\beta_*}}\) for \(n \ge 2\), and \(\rho(0) = \rho(1) = 1\). \end{cor}

A result similar to Corollary \ref{multicor} was recently established in \cite{leppanen2016} by using a more direct approach. There, the result was applied to show an almost sure ergodic theorem in the setting of quasistatic dynamical systems \cite{dobbs2016,stenlund2016}. Regarding constants, Corollary \ref{multicor} is in fact a slight improvement of the bound obtained in \cite{leppanen2016}; see Theorem 4.1 of that paper. The proof of that theorem was based on the observation that if \(h\) and \(H\) are as in Corollary \ref{multicor}, then there exist functions \(g_i \in \cC_*\) and numbers $\sigma_i \in \{-1,1\}$, such that
\begin{align*}
\widetilde{\cL}_{n_l}(Hh) = \sum_{i=1}^{2^{l+1}} \sigma_i g_i,
\end{align*}
where \(\cL_{\alpha}\) is the transfer operator associated to \(T_{\alpha}\),
\begin{align*}
\cL_{\alpha}f(x) = \sum_{y \in T_{\alpha}^{-1}x} \frac{f(y)}{T_{\alpha}'(y)},
\end{align*}
and we have denoted \(\widetilde{\cL}_{n_l} = \cL_{\alpha_{n_l}} \cdots  \cL_{\alpha_1}\). This decomposition enables one to use the known decay result of \cite{aimino2015} (Fact \ref{aimino} below) that applies to functions in the convex cone \(\cC_*\), which then leads to a multicorrelation bound similar to \eqref{eq:cor}. However, this approach produces a coefficient depending on \(l\) to the final estimate. While this does not affect the rate of decay, it is nevertheless an unnecessary dependence which is absent in the bound \eqref{eq:cor}. We also remark that Corollary \ref{multicor} applies to the class of Lipschitz continuous observables, unlike the correlation decay results of \cite{liverani1999,leppanen2016,aimino2015} where the smaller class of \(C^1\)-observables was considered.

\subsection{Results on normal approximation} The first normal approximation result of this paper follows from the main result of \cite{hella2016} (Theorem \ref{thm:stein} below) together with Theorem \ref{corr}.
 
\begin{thm}\label{stein_int_intro} Assume that \(\beta_* < 1/3\). Let \(f: \, [0,1] \to \bR^d\) be a Lipschitz continuous function with \(\hat{\mu}_{\beta_*}(f) = 0\), such that $f$ is not a coboundary in any direction\footnote{Given a unit vector $v\in\bR^d$, we say that $f$ is a coboundary in the direction~$v$ if there exists a function $g_v:[0,1]\to\bR$ in $L^2(\hat{\mu}_{\beta_*})$ such that \(v\cdot f = g_v - g_v\circ T_{\beta_*}\).}. Let \(h: \, \bR^d \to \bR\) be three times differentiable with \(\Vert D^k h \Vert_{\infty} < \infty\) for \(1 \le k \le 3\). 
Then, there is a positive-definite matrix \(\Sigma \in \bR^{d \times d}\), such that for any for \(N \ge 2\),
\begin{align*}
&\left|\hat{\mu}_{\beta_*}\left[h\left(\frac{1}{\sqrt{N}} \sum_{k=0}^{N-1} f \circ T^k_{\beta_*} \right)\right] - \Phi_{\Sigma}(h)\right|  \\
&\lesssim d^{3}\max\{1,\Vert f \Vert_{\textnormal{Lip}}^3\}(\Vert \nabla h \Vert_{\infty} + \Vert D^2 h \Vert_{\infty} + \Vert D^3 h \Vert_{\infty} ) N^{\beta_* - \frac12}(\log N)^{\frac{1}{\beta_*}} .
\end{align*}
Here \(\Phi_{\Sigma}(h)\) denotes the expectation of \(h\) with respect to the $d$-dimensional centered normal distribution $\cN(0,\Sigma)$ with covariance matrix~$\Sigma$.
\end{thm}

The second result in this vein will be established by invoking the main result of \cite{pene2005} which is again applicable by virtue of Theorem \ref{corr}. 
 
\begin{thm}\label{rio_intro} Assume that \(\beta_* < 1/3\). Let \(f: \, [0,1] \to \bR^d\) be a Lipschitz continuous function with \(\hat{\mu}_{\beta_*}(f) = 0\), such that $f$ is not a coboundary in any direction. Then, there is a positive-definite matrix \(\Sigma \in \bR^{d \times d}\), and a constant \(B(f,d,\beta_*) = B > 0\), such that for any Lipschitz continuous function \(h: \, \bR^d \to \bR\) and \(N \ge 1\),
\begin{align*}
\left|\hat{\mu}_{\beta_*}\left[h\left(\frac{1}{\sqrt{N}} \sum_{k=0}^{N-1} f \circ T^k_{\beta_*} \right)\right] - \Phi_{\Sigma}(h)\right| \le B N^{-\frac12}B\textnormal{Lip}(h).
\end{align*}
\end{thm}

The latter result gives a better rate of convergence than the former result for test functions \(h\) that are only assumed to be Lipschitz continuous. The expense is that we lose control on the dependence of the constant \(B\) with respect to \(f\) and \(d\). Both results are proved in Section \ref{multi}, where we also further discuss their differences and specify the above formulations.

\subsection*{Acknowledgements}
I would like to thank my Ph.D. advisor Mikko Stenlund for introducing to me the topic of Pomeau-Manneville maps, and for the invaluable suggestions he gave during the preparation of this paper. I gratefully acknowledge the Jane and Aatos Erkko Foundation, and the Emil Aaltosen S\"a\"ati\"o for their financial support.

\section{Preliminaries} 

The general statistical properties of the map \(T_{\alpha}\) were established by C. Liverani, B. Saussol and S. Vaienti in \cite{liverani1999}. There the authors devised a method based on a stochastic approximation of the deterministic map \(T_{\alpha}\), which enabled them to establish a polynomial rate of correlation decay. In a more recent study \cite{aimino2015}, Aimino et al. generalized the method of \cite{liverani1999} to sequences of maps \((T_{\alpha_n})_{n \ge 1}\) and showed that in this setting the polynomial correlation decay rate still applies.\footnote{Strictly speaking, the authors of \cite{aimino2015} considered a slightly modified version of the map \(T_{\alpha}\), but they pointed out that their results hold for more general maps and in particular for the map \(T_{\alpha}\). See \cite{nicol2016, aimino2015} for details.} In this section, we give a partial review of some of the developments in these two papers, emphasizing results that will be relevant for us in the subsequent sections. We also touch upon some results of our earlier paper \cite{leppanen2016}.

\indent Given \(\alpha \in (0,1)\), recall that \(\cL_{\alpha}\) denotes the transfer operator associated to \(T_{\alpha}\). We define the convex cone \(\cC_*(\beta)\) as in the subsection \ref{intermaps}. This cone is increasing,
\begin{align*}
\alpha \le \beta \:\:\Rightarrow\:\: \cC_*(\alpha) \subset \cC_*(\beta),
\end{align*}
and invariant under transfer operators in the following sense:

\begin{fact}\label{inv} If \(0 < \alpha \le \beta\), then \(\cL_{\alpha}(\cC_*(\beta)) \subset \cC_*(\beta) \). \end{fact}
 For a proof of Fact~\ref{inv}, see~\cite{liverani1999} for the original case of a single parameter and~\cite{aimino2015} for the above adaptation to a range of parameters.

Throughout this section, \((T_n)_{n\ge1}\) is a fixed admissible sequence of mappings. For \(n \ge m\), we use the following notations:
\begin{align*}
T_n &= T_{\alpha_n} \hspace{3cm} \cL_n = \cL_{\alpha_n} \\
\widetilde{T}_{n,m} &= T_n \circ \cdots \circ T_m \hspace{1cm} \widetilde{\cL}_{n,m} = {\cL}_n  \cdots  {\cL}_m \\
\widetilde{T}_n &= \widetilde{T}_{n,1} \hspace{2.95cm} \widetilde{\cL}_n = \widetilde{\cL}_{n,1}
\end{align*}

Denoting \((T_{\alpha})_1 = T_{\alpha} \upharpoonright [0,1/2]\), we have by Lemma 3.2 in \cite{liverani1999} the following estimate for the length of the leftmost branch of \(T_{\alpha}^n\):
\begin{align*}
| (T_{\alpha})_1^{-n}(0,1)| \lesssim_{\alpha} n^{-\frac{1}{\alpha}}.
\end{align*}
Then, let \((\widetilde{T}_{n,m})_1\) denote the leftmost branch of \(\widetilde{T}_{n,m}\). Since $\alpha \le \beta_*$ implies  $T_{\beta_*} \le T_{\alpha}$, it follows that
\begin{align}\label{branchlgth}
|(\widetilde{T}_{n,m})_1^{-1}(0,1)| \leq |(T_{\beta_*})_1^{-(n-m+1)}(0,1)| \lesssim (n-m+1)^{-\frac{1}{\beta_*}}.
\end{align}

\subsection{Correlation decay for Lipschitz observables.}
Let \(\rho(0) = \rho(1) = 1\), and for \(n \ge 2\) define
\begin{align*}
\rho(n) = n^{-\frac{1}{\beta_*}+1}(\log n)^{\frac{1}{\beta_*}}.
\end{align*}
The following key estimate was established in \cite{aimino2015}; see also ~\cite{liverani1999} for a similar result in the case of a single map instead of a sequence of maps.
\begin{fact}\label{aimino}
Let $f,g\in\cC_*$ with $\int f\,dx = \int g\,dx$. Then, for all $n\ge 0$,
\beqn
\|\widetilde\cL_n (f-g)\|_1 \lesssim (\|f\|_1+\|g\|_1) \rho(n) 
\eeqn
\end{fact}

Once Fact \ref{aimino} is established, the following lemma can be used for passing from the class \(\cC_*\) to Lipschitz continuous functions. 

\begin{lem}\label{lip_cone}
Suppose~$A,B\ge 0$. There exist numbers $\lambda<0$, $\nu>0$ and $\delta>0$ such that
\beqn
(f + \lambda x + \nu)h + \delta \in \cC_*
\eeqn
with
\beqn
\|(f + \lambda x + \nu)h + \delta\|_1 \lesssim AB
\eeqn
for every Lipschitz continuous function \(f\) with $\|f\|_{\textnormal{Lip}}\le A$ and every $h\in\cC_*$ with $m(h)\le B$.
In particular,
$
(\lambda x + \nu)h + \delta \in \cC_*
$
with
$
\|(\lambda x + \nu)h + \delta\|_1 \lesssim AB. 
$
\end{lem}

Lemma \ref{lip_cone} is essentially taken from \cite{liverani1999} but we provide a proof below for two reasons. First, no proof was given in \cite{liverani1999}, and secondly, the authors assumed \(f\) to be \(C^1\) in \cite{liverani1999}, which is superfluous. \\

\noindent Proof of Lemma \ref{lip_cone}: Indeed, we may choose
\beqn
\lambda = -A, \quad \nu = 6A \quad \text{and} \quad \delta = 2AB\max\!\left(\frac{a}{\beta_*+1},\frac{4a}{a-1}\right),
\eeqn
where \(a = 2^{\beta_*} (2 + \beta_*)\). Let us verify that this choice works. For brevity, denote 
\begin{align*}
F = (f + \lambda x + \nu)h + \delta.
\end{align*}

\noindent\(\mathbf{1^{\circ}}\) \(F\) is decreasing and \(F \ge 0\): Since \(\lambda \le -\text{Lip}(f) \), \(F\) is decreasing. Moreover,
\begin{align*}
f(x) + \lambda x + \nu \ge -A - A + 6A = 4A \ge 0,
\end{align*} 
so that \(F \ge 0\). \\

\noindent\(\mathbf{2^{\circ}}\) \( x^{\beta_*+1}F(x)\) is increasing: Set \(g(x) = f(x) + \lambda x + \nu\). Then, as \(g \ge 0\) and \(x^{\beta_*+1}h\) is increasing, we have for all \(y \ge x\) the lower bound
\begin{align*}
&y^{\beta_*+1}F(y) - x^{\beta_*+1}F(x)\\
&= g(y)y^{\beta_*+1}h(y) - g(x)x^{\beta_*+1}h(x) + \delta (y^{\beta_* + 1} - x^{\beta_* + 1}) \\
&= g(y)(y^{\beta_*+1}h(y) - x^{\beta_*+1}h(x)) - (g(x) - g(y))x^{\beta_*+1}h(x) + \delta (y^{\beta_* + 1} - x^{\beta_* + 1}) \\
&\ge - (g(x) - g(y))x^{\beta_*+1}h(x) + \delta (y^{\beta_* + 1} - x^{\beta_* + 1}).
\end{align*}
Notice that \(\text{Lip}(g) \le 2A\). Thus, since $0\le h(x)\le ax^{-\beta_*}m(h)$ and \(m(h) \le B\),
\begin{align*}
- (g(x) - g(y))x^{\beta_*+1}h(x) \ge -2A(y-x)axB.
\end{align*}
On the other hand,
\begin{align*}
y^{\beta_* + 1} - x^{\beta_* + 1} = (\beta_*+1)\int_x^y\xi^{\beta_*}\,d\xi \ge (\beta_*+1)(y-x)x^{\beta_*} \ge (\beta_*+1)(y-x)x.
\end{align*}
Hence, we see that
\begin{align*}
y^{\beta_*+1}F(y) - x^{\beta_*+1}F(x) \ge (y-x)x [(\beta_* +1)\delta - 2ABa] \ge 0,
\end{align*}
because
\begin{align*}
\delta \ge 2AB \frac{a}{\beta_*+1}.
\end{align*}

\noindent\(\mathbf{3^{\circ}}\) \(F(x) \le ax^{-\beta_*}m(F)\):  Note that
\beqn
4A \le \nu -\|g-\nu\|_\infty \le g \le \nu +\|g-\nu\|_\infty \le 8A.
\eeqn
Thus, using this and $h(x) \le ax^{-\beta_*}m(h)$,
\begin{align*}
F(x) &= g(x)h(x) + \delta \le  8A ax^{-\beta_*}m(h) + \delta \le 2ax^{-\beta_*}m(gh) + \delta \\
&\le (2am(gh) + \delta )x^{-\beta_*}.
\end{align*}
Since
\begin{align*}
\frac{m(gh)a}{a-1} \le \frac{8Am(h)a}{a-1} \le 2AB \frac{4a}{a-1} \le \delta,
\end{align*}
it follows that
\begin{align*}
(2am(gh) + \delta )x^{-\beta_*} \le a(m(gh) + \delta )x^{-\beta_*} = ax^{-\beta_*}m(F),
\end{align*}
as wanted. \\
\indent We have verified that \(F \in \cC_*\), and Lemma \ref{lip_cone} now obtains. \(\qed\)

Lemma \ref{lip_cone} has the following corollary.

\begin{lem}\label{recursive}
Let $f_1,f_2$ be Lipschitz continuous functions, $h\in\cC_*$ and $n\ge 0$.
Then, there exist $g_1,\dots,g_4\in\cC_*$ such that
\beqn
f_2\cdot \widetilde\cL_n(f_1 h) = g_1 - g_2 + g_3 - g_4
\eeqn
and
\beqn
\|g_i\|_1 \lesssim \|f_1\|_{\textnormal{Lip}}\|f_2\|_{\textnormal{Lip}}m(h).
\eeqn
\end{lem}

We included a proof for a version of Lemma \ref{recursive} for \(C^1\)-functions \(f_1,f_2\) in our earlier paper \cite{leppanen2016}, and this proof was based on a version of Lemma \ref{lip_cone} for \(C^1\)-functions (see Lemmas 3.1 and 3.3 in that paper). Since the exact same argument works to show Lemma \ref{recursive} (one has to only replace \(C^1\)-norms by Lipschitz-norms in the proof), we choose to omit the proof of Lemma \ref{recursive} here. This lemma can then be used to show the following version of Fact \ref{aimino} for Lipschitz continuous functions; trivial modifications to the proof of Theorem 4.1 in \cite{leppanen2016} suffice for this.

\begin{lem} Let $f,g$ be Lipschitz continuous functions with $m(f) = m(g)$. Then, for all $n\ge 0$,
\beqn
\|\widetilde\cL_n (f-g)\|_1 \lesssim (\|f\|_{\textnormal{Lip}}+\|g\|_{\textnormal{Lip}}) \rho(n).
\eeqn

\end{lem}

\subsection{The perturbed transfer operator}

\indent For \(\ve > 0\) and \(x \in [0,1)\), we denote
\begin{align*}
B_{\ve}(x) = \{ y \in [0,1) \: : \: d(x,y) \le \ve\},
\end{align*}
where \(d\) is the natural metric
\begin{align*}
d(x,y) = \min \{|x-y|, 1-|x-y|\}, \hspace{0.5cm} x,y \in [0,1).
\end{align*}
 Then, following \cite{liverani1999,aimino2015}, we introduce the family of operators
\begin{align*}
\bL_{\ve,m} = \cL_{m+n_{\ve}-1}\cdots \cL_m\bA_{\ve}, \hspace{0.5cm} m \ge 1,
\end{align*}
where \(\bA_{\ve}\) is the averaging operator
\begin{align*}
\bA_{\ve}f(x) = \frac{1}{2\ve} \int_{B_{\ve}(x)} f \, dm,
\end{align*}
and \(n_{\ve} \ge 0\) is an integer given by Fact \ref{ker} below. The operators \(\bL_{\ve,n}\) are called perturbed transfer operators. 
\begin{fact}\label{ker}  There exist \(\omega \in (0,1)\) and for all \(\ve > 0\) a number \(n_{\ve} \sim \ve^{-\beta_*}\) such that for all \(x,z \in [0,1]\) and for any \(m \ge 1\),
\begin{align*}
\frac{1}{2\ve} \widetilde{\cL}_{n_{\ve} + m -1,m} \mathbf{1}_{B_{\ve}(z)}(x) \ge \omega .
\end{align*}
\end{fact} 

Fact \ref{ker} was proved in \cite{liverani1999} in the case of a single map, and  then generalized in \cite{aimino2015} to the above form. It implies the following estimate which demonstrates how applying a random perturbation suppresses the intermittency effect; for a proof, see footnote 6 of \cite{liverani1999}.

\begin{fact}\label{exp} Let \(\omega\) be as in Fact \ref{ker}. Then, for all \(g \in L^1([0,1])\) with \(m(g) = 0\) and for all \(\ve \in (0,1)\), integers \(k,m \ge 1\),
\begin{align*}
\Vert \bL_{\ve, (k-1)n_{\ve} + m}\cdots\bL_{\ve,m}g \Vert_1 \le \exp(-\omega k)\Vert g \Vert_1.
\end{align*}
 \end{fact} 
 
\indent When operating on functions belonging to the convex cone \(\cC_*\), the \(L^1\)-distance between a usual transfer operator and its perturbed version decays polynomially in \(\ve\):

\begin{fact}\label{perturbdis} Let \(h \in \cC_*\). Then, for all integers \(m,p \ge 1\),
\begin{align*}
\Vert \widetilde{\cL}_{pn_{\ve}+m-1,m}(h) - \bL_{\ve,(p-1)n_{\ve} + m}\cdots\bL_{\ve,m}(h) \Vert_1 \lesssim p\Vert h \Vert_1\ve^{1-\beta_*}.
\end{align*}
\end{fact}

This is proved on p. 9 of \cite{aimino2015}, using a similar argument as in the proof of Lemma \ref{aux1} below.

\subsection{Densities of conditional measures.}\label{dens} Recall that \((T_n)_{n\ge1}\) is a fixed admissible sequence of maps. For each \(n \ge 1\), there is a partition \(\cP = \{ I_{\theta}^n \}_{\theta=1}^{2^n} \) of \((0,1)\) into open subintervals \(I_{\theta}^n\) such that \(\widetilde{T}_n \upharpoonright I_{\theta}^n\) maps \(I^n_{\theta}\) one-to-one and onto \((0,1)\). If \(I_1^n\) denotes the interval whose left endpoint is zero, then it follows from the definition of \(T_{\alpha}\) that
\begin{align}\label{interlength_0}
m(I_{\theta}^n) \le m(I_1^n) \hspace{0.5cm} \forall \theta \in \{1,\ldots, 2^n\}.
\end{align}
For the sake of completeness, let us verify this by induction on \(n\). \\
\indent For \(n=1\) the inequality is trivial. Suppose that \(n > 1\) and we have established \eqref{interlength_0} for \(n-1\). Let \(\theta \in \{2,\ldots, 2^n\}\). Then, we have \(m(T_1(I_{\theta}^n)) \le m(T_1(I_{1}^n))\) by the induction hypothesis (applied to the partition which corresponds to \(\widetilde{T}_{n,2}\)). We consider two cases.

\noindent\(1^{\circ}\) \(I_{\theta}^n \subset (\tfrac12, 1)\): Since \(T_1'(x) = 2 \) for all \(x\in (\tfrac12,1)\), we have \(2 m(I_{\theta}^n) = m(T_1(I_{\theta}^n)) \le m(T_1(I_{1}^n)) \le 2m(I^n_1)\). The last inequality holds because $T_1(r) = r(1 + (2r)^{\alpha_1}) \le 2r$ holds for all $r \in (0,\tfrac12)$.

\noindent\(2^{\circ}\) \(I_{\theta}^n \subset (0,\tfrac12)\): Now, \(m(T_1(I_{1}^n)) = T_1'(\xi) m(I_{1}^n) \) for some \(\xi \in I_1^n\) and \(m(T_1(I_{\theta}^n)) = T_1'(\chi) m(I_{\theta}^n) \) for some \(\chi \in I_{\theta}^n\). Since \(T_1'\) is increasing on \((0,\tfrac12)\), it follows that
\begin{align*}
m(I_{\theta}^n) \le \frac{(T_1)'(\xi) }{(T_1)'(\chi)} m(I_1^n) \le m(I_1^n).
\end{align*}
The induction is complete.

When we combine \eqref{interlength_0} with \eqref{branchlgth}, we obtain the bound 
\begin{align}\label{interlength}
m(I^n_{\theta}) \le m(I^n_1) = m((\widetilde{T}_n)^{-1}_1(0,1)) \lesssim n^{-\frac{1}{\beta_*}} \hspace{0.5cm} \forall\theta \in \{1,\ldots,2^n\}.
\end{align}

Let us fix a probability measure \(\mu\) with density \(h \in \cC_*\). We define the conditional densities
\begin{align*}
h_{\theta} = \mu(I_{\theta}^n)^{-1} \mathbf{1}_{I_{\theta}^n} h, \hspace{0.5cm} \theta \in \{1,\ldots,2^n\},
\end{align*}
and also denote \( \widetilde{h}_{\theta} = \cL_{n} \cdots \cL_{1}(h_{\theta}) = \widetilde{\cL}_{n}(h_{\theta})\). Notice that
\begin{align}\label{densform}
\widetilde{h}_{\theta}(x) = \frac{1}{\mu(I^n_{\theta})}
 \frac{h((\widetilde{T}_{n}\upharpoonright I_{\theta}^n)^{-1}x)}{(\widetilde{T}_{n})'((\widetilde{T}_{n}\upharpoonright I_{\theta}^n)^{-1}x)} .
\end{align}
Since the maps \((\widetilde{T}_{n}\upharpoonright I_{\theta}^n)^{-1} \) and \( \widetilde{T}_{n}' \circ (\widetilde{T}_{n}\upharpoonright I_{\theta}^n)^{-1}\) are increasing and \(h\) is decreasing, we see from the form \eqref{densform} that \(\widetilde{h}_{\theta}\) is decreasing.\\
\indent Fact \ref{aimino} tells us how \(\Vert \widetilde{\cL}_{m+n}(h-g)\Vert_1 \) decays as \(m \to \infty\), when \(h,g\) are probability densities belonging to \(\cC_*\). We will next show that if the initial density \(h\) is replaced by a conditional density \(h_{\theta}\), we still get a good rate of decay. To this end, we observe that the conditional densities satisfy a bound similar to that of Fact \ref{perturbdis}:

\begin{lem}\label{aux1} Let \(h \in \cC_*\) be any density. Then, for any \(m,p \ge 1\),
\begin{align*}
&\Vert \widetilde{\cL}_{pn_{\ve}+n+m,n+m+1}\widetilde{\cL}_{n+m,n+1}(\widetilde{h}_{\theta}) - \bL_{\ve,(p-1)n_{\ve} + n + m+1}\cdots\bL_{\ve,n+m+1}\widetilde{\cL}_{n+m,n+1}(\widetilde{h}_{\theta}) \Vert_1 
\\ &\lesssim \sum_{i=1}^{p}\frac{\mu(I^{n}_{\theta} \cap \widetilde{T}_{(i-1)n_{\ve} + n + m}^{-1}(0,2\ve))}{\mu(I^{n}_{\theta})}.
\end{align*} \end{lem}

\begin{myproof} The result follows by estimating as in the proof of Theorem 1.6 in \cite{aimino2015}. \\
\indent For brevity, let us denote \(A_i = \bL_{\ve, (i-1)n_{\ve} + n + m + 1} \), and \(B_i = \widetilde{\cL}_{in_{\ve} + n + m ,(i-1)n_{\ve} + n + m + 1} \), where \(i \in \{ 1,\ldots,p \}\). Then,
\begin{align*}
&\bL_{\ve,(p-1)n_{\ve} + n + m+1}\cdots\bL_{\ve,n+m+1} - \widetilde{\cL}_{pn_{\ve}+n+m,n+m+1} = A_p \cdots A_1 - B_p \cdots B_1 \\
&= \sum_{i=1}^p \left[\prod_{k=-1}^{p-i-1} A_{p-k} \right](A_i - B_i)\left[ \prod_{j=0}^{i-1} B_{i-1-j} \right],
\end{align*}
where \(A_{p+1} = 1 = B_0\). \\
\indent Fix \(i \in \{1,\ldots, p\}\) and denote \(\varphi = B_{i-1}\cdots B_{0}\widetilde{\cL}_{n+m,n+1}(\widetilde{h}_{\theta}) \). Then, \(\varphi\) is decreasing. The operators \(A_k\) are \(L^1\)-contractions, so that
\begin{align*}
\Vert A_{p+1} \cdots A_{i +1}(A_{i} - B_{i})\varphi \Vert_1 
\le \Vert (A_{i} - B_{i}) \varphi \Vert_1
\end{align*}
Moreover, since \(\varphi\) is decreasing,
\begin{align*}
&\Vert (A_{i} - B_{i}) \varphi \Vert_1 
\le \Vert \bA_{\ve}\varphi - \varphi \Vert_1 \\
&= \int_{0}^1 \left| \frac{1}{2\ve}\int_{B_{\ve}(x)} \varphi(y) \, dy - \varphi(x)  \right| \,dx \\
&\le \frac{1}{2\ve} \int_{\ve}^{1-\ve} \int_{B_{\ve}(x)} |\varphi(y) - \varphi(x)| \, dy \, dx +\int_{B_{\ve}(0)} \left| \frac{1}{2\ve} \int_{B_{\ve}(x)} \varphi(y) \, dy - \varphi(x) \right| \, dx \\
&\le \frac{1}{2\ve} \int_{\ve}^{1-\ve} \int_{x-\ve}^x \varphi(y) - \varphi(y+\ve) \, dy \, dx + 2\int_{B_{2\ve}(0)} \varphi(x) \, dx.
\end{align*}
Using Fubini's Theorem, we see that
\begin{align*}
&\frac{1}{2\ve} \int_{\ve}^{1-\ve} \int_{x-\ve}^x \varphi(y) - \varphi(y+\ve) \, dy \, dx \\
&= \frac{1}{2\ve} \int_{0}^{1-\ve} \int_{\ve}^{1-\ve} \mathbf{1}_{(y,y+\ve)}(x)( \varphi(y) - \varphi(y+\ve)) \, dx \, dy \\
&= \frac{1}{2\ve} \int_{0}^{1-\ve} m((\ve,1-\ve) \cap (y, y+\ve))( \varphi(y) - \varphi(y+\ve))  \, dy \\
& \le \frac{1}{2} \int_0^{1-\ve} \varphi(y) - \varphi(y+\ve) \, dy \le \frac12 \int_0^{\ve} \varphi(y) \, dy.
\end{align*}
By duality,
\begin{align*}
& \frac{1}{2} \int_0^{\ve} \varphi(y) \, dy = \frac{1}{2} \int_{0}^1 \widetilde{h}_{\theta} \cdot \mathbf{1}_{(0,\ve)} \circ \widetilde{T}_{(i-1)n_{\ve} + n + m,n+1} \, dm 
= \frac{1}{2} \int_{0}^1 h_{\theta} \cdot \mathbf{1}_{(0,\ve)} \circ \widetilde{T}_{(i-1)n_{\ve} + n + m} \, dm \\
&= \frac{1}{2}\frac{1}{\mu(I^{n}_{\theta})} \int_{0}^1 \mathbf{1}_{I^{n}_{\theta}} \cdot h \cdot \mathbf{1}_{(0,\ve)} \circ \widetilde{T}_{(i-1)n_{\ve} + n + m} \, dm 
= \frac{1}{2} \frac{\mu(I^{n}_{\theta} \cap \widetilde{T}_{(i-1)n_{\ve} + n + m}^{-1}(0,\ve))}{\mu(I^{n}_{\theta})}.
\end{align*}
Similarly,
\begin{align*}
2\int_{B_{2\ve}(0)} \varphi(x) \, dx \le 4\int_{0}^{2\ve} \varphi(x) \, dx 
\le  4\frac{\mu(I^{n}_{\theta} \cap \widetilde{T}_{(i-1)n_{\ve} + n + m}^{-1}(0,2\ve))}{\mu(I^{n}_{\theta})}.
\end{align*}
The bound of Lemma \ref{aux1} now follows.
\end{myproof}

Combining Lemma \ref{aux1} with the earlier estimates now yields the following:

\begin{lem}\label{condec} Let \(h,g\) be densities in \(\cC_*\). Let \(m \ge 1\), and write \(m = pn_{\ve} + l\), where \(l < n_{\ve}\). Then, for all \(\ve > 0\), \(\theta \in \{1,\ldots,2^n\}\),
\begin{align*}
&\Vert \widetilde{\cL}_{n + m,n+1}(\widetilde{h}_{\theta} - \widetilde{\cL}_{n}g) \Vert_1 \\
&\lesssim p\ve^{1-\beta_*} 
+ \sum_{i=1}^{p}\frac{\mu(I^{n}_{\theta} \cap \widetilde{T}_{(i-1)n_{\ve} + n + l}^{-1}(0,2\ve))}{\mu(I^{n}_{\theta})}
+ \exp(-\omega p).
\end{align*}
 \end{lem}

\begin{myproof} First of all, 
\begin{align*}
&\Vert \widetilde{\cL}_{n + m,n+1}(\widetilde{h}_{\theta} - \widetilde{\cL}_ng) \Vert_1 \\
&\le \Vert \widetilde{\cL}_{n + m}(g) - \bL_{\ve,(p-1)n_{\ve} + n + l +1}\cdots\bL_{\ve,n+l+1}\widetilde{\cL}_{n+l}(g) \Vert_1 \\
&+ \Vert \widetilde{\cL}_{n + m,n+1}(\widetilde{h}_{\theta}) - \bL_{\ve,(p-1)n_{\ve} + n + l +1}\cdots\bL_{\ve,n+l+1}\widetilde{\cL}_{n+l,n+1}(\widetilde{h}_{\theta}) \Vert_1 \\
&+ \Vert \bL_{\ve,(p-1)n_{\ve} + n+ l +1}\cdots\bL_{\ve,n+l+1}\widetilde{\cL}_{n+l,n+1}(\widetilde{h}_{\theta}- \widetilde{\cL}_{n}g) \Vert_1.
\end{align*}
As \(g\in \cC_*\), we have by Fact \ref{perturbdis} the estimate
\begin{align*}
&\Vert \widetilde{\cL}_{n + m}(g) - \bL_{\ve,(p-1)n_{\ve} + n + l +1}\cdots\bL_{\ve,n+l+1}\widetilde{\cL}_{n+l}(g) \Vert_1 \\
&\lesssim p\Vert \widetilde{\cL}_{n+l} g \Vert_1\ve^{1-\beta_*} \lesssim p\ve^{1-\beta_*}, \\
\end{align*}
while Lemma \ref{aux1} yields
\begin{align*}
&\Vert \widetilde{\cL}_{n + m,n+1}(\widetilde{h}_{\theta}) - \bL_{\ve,(p-1)n_{\ve} + n + l +1}\cdots\bL_{\ve,n+l+1}\widetilde{\cL}_{n+l,n+1}(\widetilde{h}_{\theta}) \Vert_1 \\
&=\Vert \widetilde{\cL}_{pn_{\ve} + n + l,n+l +1}\widetilde{\cL}_{n+l,n+1}(\widetilde{h}_{\theta}) - \bL_{\ve,(p-1)n_{\ve} + n + l +1}\cdots\bL_{\ve,n+l+1}\widetilde{\cL}_{n+l,n+1}(\widetilde{h}_{\theta}) \Vert_1 \\
 &\lesssim \sum_{i=1}^{p}\frac{\mu(I^{n}_{\theta} \cap \widetilde{T}_{(i-1)n_{\ve} + n + l}^{-1}(0,2\ve))}{\mu(I^{n}_{\theta})}.
\end{align*}
Finally, by Fact \ref{exp},
\begin{align*}
\Vert \bL_{\ve,n+(p-1)n_{\ve} + l +1}\cdots\bL_{\ve,n+l+1}\widetilde{\cL}_{n+l}(h_{\theta}-g) \Vert_1
 \le  2\exp(-\omega p).
\end{align*}
Lemma \ref{condec} now follows by putting together these three estimates.
\end{myproof}

\section{Proof of theorem \ref{corr}}

 We shall first establish the special case where \(p=1\). The general case then follows by induction on \(p\). So assume that \(p=1\) and for simplicity denote \(l = l_1\). Then, the function \(H\) in Theorem \ref{corr} becomes
\begin{align*}
H(x,y) = F(x,\widetilde{T}_{n_1}(x),\ldots, \widetilde{T}_{n_{l}}(x), \widetilde{T}_{n_{{l}+1}}(y),\ldots,\widetilde{T}_{n_{k}}(y)),
\end{align*}
where \(F\) is Lipschitz in the first \(l +1\) coordinates, and bounded in the last \(k - l\) coordinates. Without harming generality, we may assume that \(n_{l+1}-n_l \ge 2\). Denote \(n_* = n_l + \round{(n_{l+1} - n_l)/2}\)\footnote{We denote by \(\round{x}\) the greatest non-negative integer \(n\) with \(n \le x\).}, and let \( \{I_{\theta}^{n_*}\}_{\theta=1}^{2^{n_*}} \) be the partition of \((0,1)\) described at the beginning of subsection \ref{dens}. To keep notations simple, we abbreviate \(I_{\theta} = I_{\theta}^{n_*} \). Recall that the conditional densities are defined by
\begin{align*}
h_{\theta} = \mu(I_{\theta})^{-1} \mathbf{1}_{I_{\theta}} h, \hspace{0.5cm} \theta \in \{1,\ldots,2^{n_*}\},
\end{align*}
and \( \widetilde{h}_{\theta} = \widetilde{\cL}_{n_*}(h_{\theta})\). We split the domain of integration \([0,1]\) into the subintervals \(I_{\theta}\), and write
\begin{align*}
&\int H(x,x) \, d\mu(x) - \iint H(x,y) \, d\mu(x) \, d\mu_1(y)  \\
&= \int \left( H(x,x) - \int H(x,y) \, d\mu_1(y) \right) \, d\mu(x) \\
&= \sum_{\theta =1}^{2^{n_*}} \int_{I_{\theta}} \left( H(x,x) - \int H(x,y) \, d\mu_1(y) \right) \, d\mu(x).
\end{align*}
Notice that the functions \(H(\cdot,x) - \int H(\cdot,y) \, d\mu_1(y)\) are nearly constant on the intervals \(I_{\theta}\):

\begin{claim}\label{cl1} If \(a,b \in I_{\theta}\), then
\begin{align*}
&\left|H(a,x) - \int H(a,y) \, d\mu_1(y) - \left(H(b,x) - \int H(b,y) \, d\mu_1(y)\right) \right| \\
&\lesssim L (n_{l+1} - n_l)^{-\frac{1}{\beta_*}+ 1},
\end{align*}
where \(L = \max_{1 \le\alpha \le l+1} \textnormal{Lip}(F;\alpha)\).
 \end{claim}
\noindent\emph{Proof of Claim \ref{cl1}:} By induction on \(l\), we see that 
\begin{align*}
|F(a_1,\ldots, a_{l+1}, x,\ldots, x) - F(b_1,\ldots, b_{l+1}, x,\ldots, x)| \le L\sum_{i =1}^{l+1} |a_i - b_i|.
\end{align*}
Consequently,
\begin{align*}
&|H(a,x) - H(b,x) | \le L \sum_{i =0}^{l} |\widetilde{T}_{n_i}(a) - \widetilde{T}_{n_i}(b)| \leq L \sum_{i=0}^{l} m(\widetilde{T}_{n_i}I_{\theta}). 
\end{align*}
Since \(\widetilde{T}_{n_*,n_i+1}\) maps the interval \( \widetilde{T}_{n_i}I_{\theta}\) bijectively onto \((0,1)\), we have \(m(\widetilde{T}_{n_i}I_{\theta}) \lesssim (n_* - n_i)^{-\frac{1}{\beta_*}}\) by estimate \eqref{interlength}. Hence,
\begin{align*}
|H(a,x) - H(b,x) | 
&\lesssim L \sum_{i=0}^l (n_* - n_i)^{-\frac{1}{\beta_*}} 
\lesssim L \sum_{i=\round{(n_{l+1}-n_l)/2}}^{n_*} i^{-\frac{1}{\beta_*}} \\
& \lesssim  L \left( \round{\frac{n_{l+1} - n_l}{2}}^{-\frac{1}{\beta_*}} + \int^{\infty}_{\round{(n_{l+1}-n_l)/2}} t^{-\frac{1}{\beta_*}} \, dt \right)\\
&\lesssim L (n_{l+1} - n_l)^{-\frac{1}{\beta_*} + 1}.
\end{align*}
Claim \ref{cl1} now obtains. \(\qed\) \\

By Claim \ref{cl1}, fixing any \(c_{\theta} \in I_{\theta}\) for \(\theta \in \{1,\ldots,2^{n_*}\}\), we can now estimate
\begin{align}
&\left| \sum_{\theta =1}^{2^{n_*}} \int_{I_{\theta}} \left( H(x,x) - \int H(x,y) \, d\mu_1(y) \right) \, d\mu(x) \right| \notag \\
&\lesssim \left| \sum_{\theta =1}^{2^{n_*}} \int_{I_{\theta}} \left( H(c_{\theta},x) - \int H(c_{\theta},y) \, d\mu_1(y) \right) \, d\mu(x) \right| \label{main_eq1} \\
&+ L (n_{l+1} - n_l)^{-\frac{1}{\beta_*}+ 1}. \notag
\end{align}
Let \(h \in \cC_*\) denote the density of \(\mu\), and let \(h_1 \in \cC_*\) denote the density of \(\mu_1\). Moreover, let
 \(\widetilde{H}(c_{\theta},x)\) be the mapping that satisfies \(\widetilde{H}(c_{\theta}, \widetilde{T}_{n_{l+1}}(x)) = H(c_{\theta}, x) \). Then, we can rewrite the expression \eqref{main_eq1} as 
\begin{align}
&\left| \sum_{\theta =1}^{2^{n_*}} \mu(I_{\theta})\int \left( \widetilde{H}(c_{\theta},x) - \int H(c_{\theta},y) \, d\mu_1(y) \right) \, \frac{\widetilde{\cL}_{n_{l+1}}[1_{I_{\theta}}h]x}{\mu(I_{\theta})} \,dx \right| 
 \notag \\
 &= \left| \sum_{\theta =1}^{2^{n_*}} \mu(I_{\theta})\int \left( \widetilde{H}(c_{\theta},x) - \int H(c_{\theta},y) \, d\mu_1(y) \right) \, \widetilde{\cL}_{n_{l+1},n_*+1}(\widetilde{h}_{\theta})x \,dx \right|  \notag \\
&= \left| \sum_{\theta =1}^{2^{n_*}} \mu(I_{\theta})\int \left( \widetilde{H}(c_{\theta},x) - \int H(c_{\theta},y) \, d\mu_1(y) \right) \, \widetilde{\cL}_{n_{l+1},n_*+1}(\widetilde{h}_{\theta} - \widetilde{\cL}_{n_*}h_1)x \,dx \right|. \label{main_eq2}
\end{align}
We apply Lemma \ref{condec} to control the remaining sum. To this end, let \(u \ge 0\) and \(0 \le v < n_{\ve}\) be integers such that \(n_{l+1} - n_* = un_{\ve} + v\). Then, by Lemma \ref{condec},
\begin{align*}
\Vert \widetilde{\cL}_{n_{l+1},n_*+1}(\widetilde{h}_{\theta} - \widetilde{\cL}_{n_*}h_1) \Vert_1
\lesssim u\ve^{1-\beta_*} + \sum_{\gamma=1}^{u}\frac{\mu(I_{\theta} \cap \widetilde{T}_{n_* + v + (\gamma-1) n_{\ve}}^{-1}(0,2\ve))}{\mu(I_{\theta})} + \exp(-\omega u).
\end{align*}
Thus, it follows that
\begin{align*}
\eqref{main_eq2}
&\le \sum_{\theta = 1}^{2^{n_*}} \notag \mu(I_{\theta})2\Vert F \Vert_{\infty} \Vert  \widetilde{\cL}_{n_{l+1},n_*+1}(\widetilde{h}_{\theta} - \widetilde{\cL}_{n_*}h_1) \Vert_{1}  \\
&\lesssim \sum_{\theta = 1}^{2^{n_*}} \notag \mu(I_{\theta})2\Vert F \Vert_{\infty} \left[ u\ve^{1-\beta_*} + \sum_{\gamma=1}^{u}\frac{\mu(I_{\theta} \cap \widetilde{T}_{n_* + v + (\gamma-1) n_{\ve}}^{-1}(0,2\ve))}{\mu(I_{\theta})} + \exp(-\omega u) \right] \\
& \lesssim \Vert F \Vert_{\infty}(u\ve^{1-\beta_*} + \exp(-\omega u)) + \Vert F \Vert_{\infty}\sum_{\theta = 1}^{2^{n_*}}\sum_{\gamma=1}^{u}\mu(I_{\theta} \cap \widetilde{T}_{n_* + v + (\gamma-1) n_{\ve}}^{-1}(0,2\ve)).
\end{align*}
By changing the order of summation in the remaining double sum, and then using the invariance property of the cone \(\cC_*\) (i.e. Fact \ref{inv}), we see that
\begin{align*}
&\sum_{\theta = 1}^{2^{n_*}}\sum_{\gamma=1}^{u}\mu(I_{\theta} \cap \widetilde{T}_{n_* + v + (\gamma-1) n_{\ve}}^{-1}(0,2\ve)) 
= \sum_{\gamma=1}^{u}\mu( \widetilde{T}_{n_* + v + (\gamma-1) n_{\ve}}^{-1}(0,2\ve)) \\
&= \sum_{\gamma=1}^{u} \int 1_{(0,2\ve)} \circ \widetilde{T}_{n_* + v + (\gamma-1) n_{\ve}} h \, dm = \sum_{\gamma = 1}^{u} \int_0^{2\ve} \widetilde{\cL}_{n_* + v + (\gamma-1) n_{\ve}}(h) \, dm
\lesssim u \int_0^{2\ve} x^{-\beta_*} \, dx \\
&\lesssim u \ve^{1-\beta_*}.
\end{align*}
Hence, we arrive at the bound
\begin{align*}
\eqref{main_eq2} \lesssim \Vert F \Vert_{\infty}(u \ve^{1-\beta_*} + \exp(-\omega u)).
\end{align*}
Finally, since \( n_{\ve} \sim \ve^{-\beta_*}\) (see Fact \ref{ker}), 
\begin{align*}
 \eqref{main_eq2} &\lesssim \Vert F \Vert_{\infty}(u \ve^{1-\beta_*} + \exp(-\omega u))  \\
&\lesssim \Vert F \Vert_{\infty} \left( (n_{l+1}-n_*) \ve + \exp\left(-\omega \ve^{\beta_*} (n_{l+1}-n_*) \right) \right)  \\
&\lesssim \Vert F \Vert_{\infty} \left( \frac{n_{l+1}-n_l}{2} \ve + \exp\left(-\omega \ve^{\beta_*} \frac{n_{l+1}-n_l}{2} \right) \right)  \\
&\lesssim \Vert F \Vert_{\infty} \rho(n_{l+1} - n_l),
\end{align*}
when we choose \(\ve = (n_{l+1}-n_l)^{-\frac{1}{\beta_*}} \left( \log (n_{l+1} - n_l)^{\kappa (\frac{1}{\beta_*} -1)}\right)^{\frac{1}{\beta_*}} \) for appropriate \(\kappa > 0\). This completes the proof of the case \(p = 1\). 

Then, assume that we have already obtained the bound of Theorem \ref{corr}  for \(1 \le l_1 \le \ldots \le l_{p-1} < k\), and suppose that \(l_{p-1} < l_{p} < k\). Let \(F : \, [0,1]^{k+1} \to \bR\) be a bounded function that is Lipschitz continuous in the coordinates \(x_{\alpha}\) for \(1 \le \alpha \le l_p +1\),  and recall that \(H(x_0,\ldots, x_{p})\) denotes the function
\begin{align*}
F(x_0,\widetilde{T}_{n_1}(x_0),\ldots, \widetilde{T}_{n_{l_1}}(x_0), \widetilde{T}_{n_{{l_1}+1}}(x_1),\ldots,\widetilde{T}_{n_{{l_2}}}(x_1), \ldots, \widetilde{T}_{n_{l_p+1}}(x_{p}),\ldots \widetilde{T}_{n_k}(x_{p})).
\end{align*}
From the case \(p=1\), we know that
\begin{align}
&\left| \int H(x,\ldots,x) \, d\mu(x)  - \iint H(x,\ldots x,x_{p}) \, d\mu (x) \, d\mu_p (x_{p}) \right| \notag \\
&\lesssim 
(\Vert F \Vert_{\infty} + L) \rho(n_{l_{p}+1}-n_{l_p}), \label{est1}
\end{align}
where \(L = \max_{1 \le\alpha \le l_p+1} \textnormal{Lip}(F;\alpha)\). \\
\indent Then, for each \(x_{p} \in [0,1]\), we can apply the induction hypothesis to the function \((y_0,\ldots,y_k) \mapsto F(y_0,\ldots,y_{l_p},\widetilde{T}_{n_{l_p+1}}(x_{p}), \ldots,\widetilde{T}_{n_k}(x_{p}))\), which is bounded and Lipschitz continuous in all of its coordinates. This yields the estimate
\begin{align}
&\left| \int H(x,\ldots x,x_{p}) \, d\mu (x) - \idotsint H(x_0,\ldots, x_{p-1},x_p) \, d\mu(x_0) \, d\mu_1(x_1) \ldots d\mu_{p-1}(x_{p-1}) \right| \notag \\
&\lesssim ( \Vert F \Vert_{\infty} + \max_{1 \le\alpha \le l_{p-1}+1} \textnormal{Lip}(F;\alpha))\sum_{i=1}^{p-1}\rho(n_{l_{i}+1}-n_{l_i}), \label{est2}
\end{align}
for all \(x_p \in [0,1]\). Now, to obtain the desired bound, it suffices to apply the triangle inequality in combination with \eqref{est1} and \eqref{est2}. \\
\indent This completes the proof of Theorem \ref{corr}.
\(\qed\)
\section{Multivariate normal approximation}\label{multi}

In this section, we briefly review the general theories of \cite{hella2016, pene2005} and show how they can be applied together with Theorem \ref{corr} to establish the normal approximation results \ref{stein_int_intro} and \ref{rio_intro}. We begin by introducing some further notation. Let us remark that unlike in the previous sections, we will from now on consider a single map only instead of a sequence of maps. 

 Let \(T:\, X \to X\) be a measure preserving transformation on a probability space \((X, \cB, \mu)\), and let \(f : \, X \to \bR^d\) be a given function. Then, for \(k \in \bN_0 = \{0,1,\ldots\}\), \(T^k\) denotes the \(k\)-fold composite of \(T\), where \(T^0 = \text{id}_X\), and \(f^k = f \circ T^k\). For \(N\in \bN_0\), we write 
\beqn
W = W(N) = \frac{1}{\sqrt{N}}\sum_{k=0}^{N-1} f^k.
\eeqn
Given $K\in \bN_0\cap [0,N-1]$, let
\beqn
[n]_K = [n]_K(N) = \{ k \in \bN_0 \cap [0,N-1] \, : \, |k-n| \le K \}
\eeqn
and
\beq\label{eq:W^n}
W^n = W^n(N,K) = W - \frac{1}{\sqrt{N}} \sum_{k\in [n]_K} f^k
\eeq
for all $n\in\bN_0\cap [0,N-1]$. In other words, $W^n$ differs from~$W$ by a time gap (within $[0,N-1]$) of radius~$K$, centered at time~$n$. 

Given a unit vector $v\in\bR^d$, we say that $f$ is a coboundary in the direction~$v$ if there exists a function $g_v:X\to\bR$ in $L^2(\mu)$ such that
\beqn
v\cdot f = g_v - g_v\circ T.
\eeqn

We denote by $\Phi_\Sigma(h)$ the expectation of a function $h:\bR^d\to\bR$ (in this section \(h\) does not stand for a density anymore) with respect to the $d$-dimensional centered normal distribution $\cN(0,\Sigma)$ with covariance matrix~$\Sigma\in\bR^{d\times d}$, i.e.,
\beqn
\Phi_\Sigma(h) = \frac{1}{\sqrt{(2\pi)^d\det\Sigma}} \int_{\bR^d} e^{-\frac12 w\cdot\Sigma^{-1}w}h(w)\, dw.
\eeqn

For a function \(B: \, \bR^d \to \bR\), we write \(D^kB\) for the $k$th derivative of \(B\), and also denote \(\nabla B = D^1 B\). We define
\begin{align*}
\Vert D^k B \Vert_{\infty} = \max \{ \|\partial_1^{t_1}\cdots\partial_d^{t_d}B_\alpha\|_{\infty}  :  t_1+\cdots+t_d = k,\, 1\le \alpha\le d'  \}.
\end{align*} 
Finally, given two vectors $v,w\in\bR^d$, we write $v\otimes w$ for the $d\times d$ matrix with entries
\beqn
(v\otimes w)_{\alpha\beta} = v_\alpha w_\beta.
\eeqn

\subsection{Stein's method.} Let \(\Sigma \in \bR^{d \times d}\) be a symmetric positive definite matrix. Stein's method of normal approximation gives a way of bounding the distance 
\begin{align*}
|\mu(h(W)) - \Phi_{\Sigma}(h)| 
\end{align*}
for a given function \(h: \bR^d \to \bR\)  by introducing the so called Stein equation
\begin{align}\label{eq:steinmv1}
\tr  \Sigma D^2A(w) - w \cdot  \nabla A(w) = h(w) - \Phi_{\Sigma}(h),
\end{align}
where \(\tr  \Sigma D^2A(w)\) denotes the trace of the matrix \(\Sigma D^2A(w)\). If we evaluate both sides of \eqref{eq:steinmv1} at the random vector \(W\) and integrate with respect to \(\mu\), we arrive at
\begin{align}\label{eq:steinmv2}
\mu[\tr  \Sigma D^2A(W) - W \cdot  \nabla A(W)] = \mu(h(W)) - \Phi_{\Sigma}(h).
\end{align}
Thus, we see that a possible strategy for bounding \(\mu(h(W)) - \Phi_{\Sigma}(h)\) proceeds by solving the differential equation \eqref{eq:steinmv1} and bounding the left-hand side of \eqref{eq:steinmv2}. 
If the transformation \(T\) exhibits sufficiently rapid decay of correlations, the latter task can be done by exploiting the auxiliary random vector \(W^n\). This approach involves Taylor expanding $\nabla A(W^n)$ at $W$, which calls for bounds on the partial derivatives of \(A\). After conditions on \(A\) have been derived, they can be translated to conditions on \(h\) by using the explicit solution to the equation \eqref{eq:steinmv1}. Indeed, \(\Vert D^k h \Vert_{\infty} < \infty\) for \(1 \le k \le 3\) implies that  \(A \in C^3(\bR^d,\bR)\) and \(\Vert D^k A \Vert_{\infty} < \infty\) for \(1 \le k \le 3\); see \cite{hella2016, gaunt2016}. These steps lead to the following result.

\begin{thm}\label{thm:stein} Let \(f: \, X \to \bR^d\) be a bounded measurable function with \(\mu(f) = 0\). Let \(h: \, \bR^d \to \bR\) be any three times differentiable function with \(\Vert D^k h \Vert_{\infty} < \infty\) for \(1 \le k \le 3\). Fix integers $N>0$ and $0\le K<N$. Suppose that the following conditions are satisfied:

\begin{itemize}
\item[\bf (A1)] There exist constants $C_2 > 0$ and $C_4 > 0$, and a non-increasing function \(\tau : \, \bN_0 \to \bR_+\) with \(\tau(0) = 1\) and \( \sum_{i=1}^{\infty} i\tau(i) < \infty\), such that
\begin{align*}
|\mu(f_{\alpha} f^k_{\beta})| & \le C_2\,\tau(k)
\\
|\mu(f_{\alpha} f_{\beta}^l f_{\gamma}^m f_{\delta}^n)| & \le C_4 \min\{\tau(l),\tau(n-m)\}
\\
|\mu(f_{\alpha} f_{\beta}^l f_{\gamma}^m f_{\delta}^n) - \mu(f_{\alpha} f_{\beta}^l)\mu(f_{\gamma}^m f_{\delta}^n)| & \le C_4\,\tau(m-l)
\end{align*}
hold whenever $k\ge 0$; $0\le l\le m \le n < N$; $\alpha,\beta,\gamma,\delta\in\{\alpha',\beta'\}$ and $\alpha',\beta'\in\{1,\dots,d\}$.

\smallskip
\item[\bf (A2)] There exists a function \(\tilde\tau : \, \bN_0 \to \bR_+\) such that
\begin{align*}
|\mu( f^n \cdot  \nabla h(v + W^n t) ) | \le \tilde\tau(K)
\end{align*}
holds for all $0\le n < N$, $0\le t\le 1$ and $v\in\bR^d$.

\smallskip
\item[\bf (A3)] $f$ is not a coboundary in any direction. 

\end{itemize}
Then
\begin{align}\label{variance}
\Sigma = \mu(f \otimes f) + \sum_{n=1}^{\infty} (\mu(f^n \otimes f) + \mu(f \otimes f^n))
\end{align}
is a well-defined, symmetric, positive-definite, $d\times d$ matrix; and
\begin{align}\label{eq:main}
|\mu(h(W)) - \Phi_{\Sigma}(h)| \le C_* \! \left(\frac {K+1}{\sqrt{N}} + \sum_{i= K+1}^{\infty}  \tau(i)\right)  + \sqrt N \tilde\tau(K),
\end{align}
where
\begin{align*}
C_* = 12d^3\max\{C_2,\sqrt{C_4}\}\left(\Vert D^2 h\Vert_{\infty} + \Vert f \Vert_{\infty} \Vert D^3 h \Vert_{\infty} \right)\sum_{i = 0}^\infty (i+1)\tau(i) 
\end{align*}
is independent of $N$ and $K$.
\end{thm}
 
Theorem \ref{thm:stein} was the main result of \cite{hella2016} where a proof can be found. The assumptions as well as the conclusion of the result concern only the given functions \(f\) and \(h\), and the fixed numbers \(N\) and \(K\). Moreover, the constant \(C_*\) in the bound is expressed entirely in terms of the quantities appearing in the assumptions. Condition (A1) requires decay of correlations of orders two and four, at a rate which has a finite first moment. Observe that for this to hold in the setting of the Pomeau-Manneville map \(T_{\beta_*}\), we need to require that \(\beta_* < 1/3\). Moreover, for the bound in~\eqref{eq:main} to be of any use, we need $K\ll \sqrt N$, and $\tilde\tau(K)$ has to be small.   \\
\indent Let us now return to the setting of Pomeau-Manneville maps, i.e. let \((T, X, \cB, \mu) = (T_{\beta_*}, [0,1], \text{Borel}([0,1]), \hat{\mu}_{\beta_*})\).

\begin{thm}\label{stein_int} Assume that \(\beta_* < 1/3\). Let \(f: \, [0,1] \to \bR^d\) be a Lipschitz continuous function with \(\hat{\mu}_{\beta_*}(f) = 0\), such that $f$ is not a coboundary in any direction. Let \(h: \, \bR^d \to \bR\) be three times differentiable with \(\Vert D^k h \Vert_{\infty} < \infty\) for \(1 \le k \le 3\). 
Then for \(N \ge 2\),
\begin{align*}
&|\hat{\mu}_{\beta_*}(h(W)) - \Phi_{\Sigma}(h)| \\
&\lesssim d^{3}\max\{1,\Vert f \Vert_{\textnormal{Lip}}^3\}(\Vert \nabla h \Vert_{\infty} + \Vert D^2 h \Vert_{\infty} + \Vert D^3 h \Vert_{\infty} ) N^{\beta_* - \frac12}(\log N)^{\frac{1}{\beta_*}},
\end{align*}
where \(\Sigma \in \bR^{d\times d}\) is the positive-definite matrix given by \eqref{variance}.
\end{thm}

\begin{myproof} It suffices to verify the assumptions (A1) and (A2) of Theorem \ref{thm:stein}. We do this by applying Theorem \ref{corr}. 

(A1): We choose \(\tau = \rho\). Then, under the standing assumption \(\beta_* < 1/3\), we have \( \sum_{i=1}^{\infty} i\rho(i) < \infty\).  If \(F: \, \bR^{n+1} \to \bR\), \(F(x_0,\ldots,x_{n}) = f_{\alpha_1}(x_1)\cdots f_{\alpha_{n}}(x_{n})\), where \(\alpha_i \in \{1,\ldots ,d\}\), then \(\Vert F \Vert_{\infty} \le \Vert f \Vert_{\infty}^n \le \Vert f \Vert_{\text{Lip}}^n \) and \( \text{Lip}(F; \alpha) \le \Vert f \Vert_{\text{Lip}}^n\) for all \(\alpha =1,\ldots, n+1\). Thus, by Theorem \ref{corr},
\begin{align*}
|\hat{\mu}_{\beta_*}(f_{\alpha} f^k_{\beta})| 
&\lesssim  \Vert f \Vert_{\text{Lip}}^2 \rho(k) \\
|\hat{\mu}_{\beta_*}(f_{\alpha} f_{\beta}^l f_{\gamma}^m f_{\delta}^n)| & \lesssim \Vert f \Vert_{\text{Lip}}^4 \min\{\rho(l),\rho(n-m)\}
\\
|\hat{\mu}_{\beta_*}(f_{\alpha} f_{\beta}^l f_{\gamma}^m f_{\delta}^n)) - \hat{\mu}_{\beta_*}(f_{\alpha} f_{\beta}^l)\hat{\mu}_{\beta_*}(f_{\gamma}^m f_{\delta}^n)| & \lesssim \Vert f \Vert_{\text{Lip}}^4\,\rho(m-l),
\end{align*} 
whenever $k\ge 0$; $0\le l\le m \le n < N$; $\alpha,\beta,\gamma,\delta\in\{1,\ldots ,d\}$ . 

 (A2): Let \(h: \, \bR^d \to \bR\) be three times differentiable with \(\Vert D^k h \Vert_{\infty} < \infty\) for \(1 \le k \le 3\). If $0\le n < N$, $0\le t\le 1$ and $v\in\bR^d$, then Theorem \ref{corr} implies the bound
\begin{align}\label{bound_a2}
|\hat{\mu}_{\beta_*}( f^n \cdot \nabla h(v + W^n t) ) | \lesssim d^{2}(\Vert \nabla h \Vert_{\infty} \Vert f \Vert_{\text{Lip}} + \Vert f \Vert_{\text{Lip}}^2 \Vert D^2h \Vert_{\infty})\rho(K).
\end{align}
Indeed, let 
\begin{align*}
F(x_0,\ldots, x_{n-K},x_n, x_{n+K}, \ldots, x_{N-1}) = f(x_n) \cdot \nabla h \left(v + \frac{1}{\sqrt{N}} \sum_{i \notin [n]_K} f(x_i)t \right).
\end{align*}
Then \( \Vert F \Vert_{\infty} \le  d\Vert f \Vert_{\infty} \Vert \nabla h \Vert_{\infty}\) and  for all \(\alpha \le N\),
\begin{align*}
\text{Lip}(F; \alpha) &\le d \Vert \nabla h \Vert_{\infty} \Vert f \Vert_{\text{Lip}} + d^{2} \Vert f \Vert_{\text{Lip}}^2 N^{-\frac12}\Vert D^2h \Vert_{\infty},
\end{align*}
so that Theorem \ref{corr} is applicable with \(F\). This yields \eqref{bound_a2}. \\
\indent With the foregoing bounds, it now follows by Theorem \ref{thm:stein} that
\begin{align*}
&|\hat{\mu}_{\beta_*}(h(W)) - \Phi_{\Sigma}(h)| \\
&\lesssim d^3 \Vert f \Vert_{\text{Lip}}^2\left(\Vert D^2 h\Vert_{\infty} + \Vert f \Vert_{\infty} \Vert D^3 h \Vert_{\infty} \right)  \! \left(\frac {K+1}{\sqrt{N}} + \sum_{i= K+1}^{\infty}  \rho(i)\right) \\ 
&+ \sqrt N d^{2}(\Vert \nabla h \Vert_{\infty} \Vert f \Vert_{\text{Lip}} + \Vert f \Vert_{\text{Lip}}^2 \Vert D^2h \Vert_{\infty})\rho(K) \\
& \lesssim d^{3}\max\{1,\Vert f \Vert_{\text{Lip}}^3\}(\Vert \nabla h \Vert_{\infty} + \Vert D^2 h \Vert_{\infty} + \Vert D^3 h \Vert_{\infty} )\left(\frac {K}{\sqrt{N}} + \sum_{i= K+1}^{\infty}  \rho(i) + \sqrt{N} \rho(K)\right).
\end{align*}
With \(K \le \sqrt{N}\), we have \(\sum_{i \ge K +1} i^{1-\frac{1}{\beta_*}} \lesssim K^{2-\frac{1}{\beta_*}} \lesssim \sqrt{N} K^{1-\frac{1}{\beta_*}} \), and thus we choose \(K = \round{N^{\beta_*}}\) so that \(\sqrt{N} K^{1-\frac{1}{\beta_*}} \approx K/\sqrt{N} \). Then,
\begin{align*}
&\frac {K}{\sqrt{N}} + \sum_{i= K+1}^{\infty}  \rho(i) + \sqrt{N} \rho(K) \\
&\lesssim N^{\beta_* - \frac12} + \sum_{i= K+1}^{\infty} i^{1 - \frac{1}{\beta_*}} (\log i)^{\frac{1}{\beta_*}} + N^{\beta_* - \frac12}(\log N)^{\frac{1}{\beta_*}} \\
&\lesssim N^{\beta_* - \frac12}(\log N)^{\frac{1}{\beta_*}}.
\end{align*}
Consequently,
\begin{align*}
&|\hat{\mu}_{\beta_*}(h(W)) - \Phi_{\Sigma}(h)| \\
&\lesssim d^{3}\max\{1,\Vert f \Vert_{\text{Lip}}^3\}(\Vert \nabla h \Vert_{\infty} + \Vert D^2 h \Vert_{\infty} + \Vert D^3 h \Vert_{\infty} ) N^{\beta_* - \frac12}(\log N)^{\frac{1}{\beta_*}}.
\end{align*}
This finishes the proof of Theorem \ref{stein_int}.
\end{myproof}

 The authors of \cite{nicol2016} recently established a univariate CLT in the setting of sequential Pomeau-Manneville  maps, i.e. when \(f \circ T_{\beta_*}^k\) is replaced by \(f \circ T_{\alpha_k} \circ \cdots \circ T_{\alpha_1}\), where \((T_{\alpha_k})_{k\ge 1}\) is a suitable admissible sequence of maps. We mention in passing that a virtue of Stein's method visible from \cite{hella2016} is its quite simple and conceptually transparent nature, which ultimately enables one to implement it in non-stationary situations. Indeed, the method can be used to prove a generalization of Theorem \ref{thm:stein} for sequences of transformations. The generalization can then be applied with Theorem \ref{corr} to obtain a version of Theorem \ref{stein_int} for admissible sequences \((T_{\alpha_k})_{k\ge 1}\) of Pomeau-Manneville maps with $\beta_* < 1/3$.

\subsection{Rio's method.} F. P\`{e}ne's adaptation of Rio's method gives conditions for estimating the rate of convergence in the multivariate CLT by the Kantorovich metric
\begin{align*}
d_{\kappa}(W,Z) = \sup \{| \mu(h(W)) - \Phi_{\Sigma}(h) | \: : \: h : \bR^d \to \bR, \: \text{Lip}(h) \le 1 \},
\end{align*}
where \(Z \sim \cN(0,\Sigma)\). Theorem \ref{rio} below was established in P\`{e}ne's paper \cite{pene2005}, and it shows that in the case of systems with appropriately decaying correlations \(d_{\kappa}(W,Z)\) decays at the optimal rate \(N^{-1/2}\). The reader should consult \cite{pene2005} for a detailed treatment of this method, as well as for more information about different quantities used to estimate the rate of convergence in the CLT. \\
\indent We work in the setting of a general measure preserving transformation \((X, \cB, \mu, T)\), and continue to use the notations  introduced at the beginning of this section. Additionally, for two functions \(f,g\) belonging to \(L^2(\mu)\), we denote \(\text{Cov}(f,g) = \mu(fg) - \mu(f)\mu(g)\), and  for \(N \ge 1\), we denote \(S_N = \sum_{k=0}^{N-1} f^k\) so that \(W = N^{-1/2}S_N\). We also let \(S_0 = 0\). 

\begin{thm}\label{rio} Let \(f: \, X \to \bR^d\) be a bounded function with \(\mu(f) = 0\). Assume that there exist \(r \in \bN_0\), \(C \ge 1\), \(M \ge \max\{ 1, \Vert f \Vert_{\infty} \}\) and a sequence of non-negative real numbers \((\varphi_{p,l})\), such that the following conditions hold:
\begin{itemize}
\item[\bf (B1)] \(|\varphi_{p,l}| \le 1 \) and \( \sum_{p=1}^{\infty} p \max_{0 \le l \le \round{p/(r+1)}} \varphi_{p,l} < \infty\). \\

\item[\bf (B2)] For any integers \(a,b,c\) satisfying \(1 \le a + b + c \le 3\), for any integers \(i,j,k,p,q,l\) with \(0 \le i \le j \le k \le k + p \le k + p + q \le k + p + l\), for any \(\alpha, \beta, \gamma \in \{1,\ldots , d\}\), and for any bounded differentiable function \(F : \, \bR^d \times ([-M,M]^d)^3 \to \bR\) with bounded gradient,
\begin{align}\label{cond:b2}
&| \textnormal{Cov}[ F(S_{i},f^i,f^j,f^k),(f^{k+p}_{\alpha})^a(f^{k+p+q}_{\beta})^b(f^{k+p+l}_{\gamma})^c] | \\
&\le C(\Vert F \Vert_{\infty} + \Vert \nabla F \Vert_{\infty})\, \varphi_{p,l}. \notag
\end{align}
\smallskip
\item[\bf (B3)] $f$ is not a coboundary in any direction.\footnote{This condition was not stipulated in the main result of \cite{pene2005}, but we have added it to ensure that the covariance \(\Sigma\) is positive definite. This is not necessary, if a more general definition of normal distribution, such as the one given in \cite{pene2005}, is used.}
\end{itemize}
Then
\begin{align}\label{variance2}
\Sigma = \mu(f \otimes f) + \sum_{n=1}^{\infty} (\mu(f^n \otimes f) + \mu(f \otimes f^n))
\end{align}
is a well-defined, symmetric, positive-definite, $d\times d$ matrix; and
 there exists \(B > 0\) such that for any Lipschitz continuous function \(h: \, \bR^d \to \bR\), \(N \ge 1\),
\begin{align*}
| \mu(h(W)) - \Phi_{\Sigma}(h) | \le \frac{B\textnormal{Lip}(h)}{\sqrt{N}}. 
\end{align*}
\end{thm}

In \cite{pene2005}, Theorem \ref{rio} was formulated and proved for an arbitrary centered stationary process, but for simplicity we have stated it here in the less general setting of measure preserving transformations. In contrast with Theorem \ref{thm:stein}, Theorem \ref{rio} establishes the optimal convergence rate \(N^{-1/2}\). However, as discussed in \cite{hella2016}, there is a notable difference in the constants of the two upper bounds. Due to the smooth metric used in Theorem \ref{thm:stein}, the upper bound is independent of the limiting covariance \(\Sigma\), while the constant \(B\) in Theorem \ref{rio} depends on \(\Sigma\).  E.g. if \(d=1\), even in the case of independent random variables the Kantorovich distance depends on how the limiting variance compares with higher absolute moments (the same phenomenon is seen in the classical Berry-Esseen theorem). We also remark that the explicit expression of the constant \(C_*\) in Theorem \ref{thm:stein} enables us to see that the upper bound of Theorem \ref{thm:stein} (and Theorem \ref{stein_int}) scales as \(d^3\) with increasing dimension of the observable \(f\). \\
\indent Condition (B2) of Theorem \ref{rio} is resemblant to condition (A2) of Theorem \ref{thm:stein}, in that both of them require the decay of a certain functional expression depending on a finite fragment of the trajectory of \(T\). This is again guaranteed by Theorem \ref{corr} in the setting \((T, X, \cB, \mu) = (T_{\beta_*}, [0,1], \text{Borel}([0,1]), \hat{\mu}_{\beta_*})\), resulting in the following estimate. 

\begin{thm}\label{rio_int} Assume that \(\beta_* < 1/3\). Let \(f: \, [0,1] \to \bR^d\) be a Lipschitz continuous function with \(\hat{\mu}_{\beta_*}(f) = 0\), such that $f$ is not a coboundary in any direction. Then, there exists a constant \(B(f,d, \beta_*) = B > 0\) such that for any Lipschitz continuous function \(h: \, \bR^d \to \bR\) and \(N \ge 1\),
\begin{align*}
&|\hat{\mu}_{\beta_*}(h(W)) - \Phi_{\Sigma}(h)| \le N^{-\frac12}B\textnormal{Lip}(h),
\end{align*}
where \(\Sigma \in \bR^{d\times d}\) is the positive-definite matrix given by \eqref{variance2}.
\end{thm}

\begin{myproof} Theorem \ref{rio_int} follows, once we verify the conditions (B1) and (B2) of Theorem \ref{rio}. We let \(\varphi_{p,l} = \rho(p)\) for all \(p,l \ge 0\), so that condition (B1) is satisfied, since we are assuming \(\beta_* < 1/3\). \\
\indent Then, let \(a,b,c\in \bN_0\), \(i,j,k,p,q,l \in \bN_0\) and \(\alpha, \beta, \gamma \in \{1,\ldots , d\}\) be as in the condition (B2). Moreover, fix a differentiable function \(F : \, \bR^d \times ([-\Vert f \Vert_{\infty} - 1, \Vert f \Vert_{\infty} + 1]^d)^3 \to \bR\) with bounded gradient. We need to show that for a suitable \(C > 0\),
\begin{align}\label{eq:last}
&| \textnormal{Cov}[ F(S_{i},f^i,f^j,f^k),(f^{k+p}_{\alpha})^a(f^{k+p+q}_{\beta})^b(f^{k+p+l}_{\gamma})^c] | \notag\\
&\le C(\Vert F \Vert_{\infty} + \Vert \nabla F \Vert_{\infty})\, \rho(p). 
\end{align}
This follows by Theorem \ref{corr}. To see this, define 
\begin{align*}
&G(x_0,\ldots, x_{i-1}, x_i, x_j, x_k, x_{k+p}, x_{k+p+q}, x_{k+p+l}) \\
&= F\left(\sum_{m=0}^{i-1} f(x_m), f(x_i), f(x_j),f(x_k)\right)f_{\alpha}^a(x_{k+p})f_{\beta}^b(x_{k+p+q})f_{\gamma}^c(x_{k+p+l}).
\end{align*}
Then \(\Vert G \Vert_{\infty} \le \Vert F \Vert_{\infty} \Vert f \Vert_{\infty}^{a+b+c}\), and  
\begin{align*}
\text{Lip}(G;\alpha) \le d\Vert f \Vert_{\infty}^{a+b+c}\Vert \nabla F \Vert_{\infty} \text{Lip}(f), \hspace{0.2cm} \alpha = 1,\ldots i+3.
\end{align*}
Thus, Theorem \ref{corr} can be applied to the function \(G\), which implies  \eqref{eq:last} for~$C = c({\beta_*}) \cdot d\max\{\Vert f \Vert_{\text{Lip}}^4,1\}$, where \(c({\beta_*})\) is a constant depending only on \(\beta_*\).
 \end{myproof}

\bibliography{stein_int}{}
\bibliographystyle{plainurl}

\end{document}